\documentclass[11pt, reqno]{amsart}
\usepackage{amssymb}
\usepackage{amscd}
\usepackage{verbatim,ifthen}
\usepackage{color}
\usepackage{latexsym}
\usepackage{tikz}
\usepackage{tikz-cd}
\usepackage{mathrsfs}
\usepackage{wrapfig}
\usetikzlibrary{shapes}
\usepackage{color}
\usetikzlibrary{arrows.meta}
\usepackage{bbm}
\usetikzlibrary{matrix}
\usetikzlibrary{calc}
\usetikzlibrary{arrows,intersections}
\usepackage{pgfplots}
\usepackage{multicol}
\usepackage{array}
\newcolumntype{M}[1]{>{\centering\arraybackslash}m{#1}}

\usepackage[colorlinks, linkcolor=black, citecolor=blue, linktocpage,urlcolor=blue]{hyperref}
\addtocontents{toc}{\setcounter{tocdepth}{1}}

\usepackage{amsmath}

\numberwithin{equation}{section}

\let\oldtocsection=\tocsection
 
\let\oldtocsubsection=\tocsubsection

\renewcommand{\tocsection}[2]{\hspace{0em}\oldtocsection{#1}{#2}}
\renewcommand{\tocsubsection}[2]{\hspace{1em}\oldtocsubsection{#1}{#2}}

\def\XXint#1#2#3{{\setbox0=\hbox{$#1{#2#3}{\int}$ }
\vcenter{\hbox{$#2#3$ }}\kern-.6\wd0}}

\usepackage{eucal}
\usepackage{calc}  
\usepackage{enumitem} 
\usepackage{tensor}
\usepackage{graphicx,wrapfig,lipsum}
\usepackage{etoolbox}
\usepackage{marginnote}
\usepackage{lipsum}
\makeatletter
\patchcmd{\@mn@margintest}{\@tempswafalse}{\@tempswatrue}{}{}
\patchcmd{\@mn@margintest}{\@tempswafalse}{\@tempswatrue}{}{}
\reversemarginpar 
\makeatother
\usepackage{scrextend}

\makeatletter
\DeclareRobustCommand\widecheck[1]{{\mathpalette\@widecheck{#1}}}
\def\@widecheck#1#2{%
    \setbox\z@\hbox{\m@th$#1#2$}%
    \setbox\tw@\hbox{\m@th$#1%
       \widehat{%
          \vrule\@width\z@\@height\ht\z@
          \vrule\@height\z@\@width\wd\z@}$}%
    \dp\tw@-\ht\z@
    \@tempdima\ht\z@ \advance\@tempdima2\ht\tw@ \divide\@tempdima\thr@@
    \setbox\tw@\hbox{%
       \raise\@tempdima\hbox{\scalebox{1}[-1]{\lower\@tempdima\box
\tw@}}}%
    {\ooalign{\box\tw@ \cr \box\z@}}}
\makeatother

\usepackage{soul}

\title{On Hermitian manifolds with vanishing curvature}
\author{Kyle Broder}
\address{The University of Queensland,  St. Lucia,  QLD 4067, Australia}
\email{k.broder@uq.edu.au}
\author{Kai Tang}
\address{College of Mathematics and Computer Science, Zhejiang Normal University, Jinhua, Zhejiang, 321004, China}
\email{kaitang001@zjnu.edu.cn}
\thanks{The first named author was partially supported by the Australian Government through the Australian Research Council's Discovery Projects funding scheme (project DP220102530) and an Australian Government Research Training Program (RTP) Scholarship. The second named author was supported by National Natural Science
Foundation of China (Grant No.12001490).}
\begin{document}

\maketitle

\begin{abstract}
We show that Hermitian metrics with vanishing holomorphic curvature on compact complex manifolds with pseudoeffective canonical bundle are conformally balanced. Pluriclosed metrics with vanishing holomorphic curvature on compact K\"ahler manifolds are shown to be K\"ahler and hence, are completely classified. We prove that Hermitian metrics with vanishing real bisectional curvature on complex manifolds in the Fujiki class $\mathcal{C}$ are K\"ahler and thus fall under the same classification. Finally, we formalize the notion of `altered' curvatures, which force distinguished metric structures when mandated to coincide with their `standard' counterparts.
\end{abstract}

\section{Introduction}
\noindent Hermitian geometry has exposed a vast, untamed wilderness of curvature conditions.  Even within K\"ahler geometry,  the relationship between some of the most classical curvatures remains far from understood (see, e.g.,  \cite{BroderRemarks}).  A foundational result in the subject is the classification of simply connected complete K\"ahler manifolds with constant holomorphic sectional curvature $\text{HSC}_{\text{g}} \equiv c$,  where \begin{eqnarray}\label{HSCExpression}
\text{HSC}_{\text{g}}(\xi) & : = & \frac{1}{| \xi |_{\text{g}}^4} \sum_{i,j,k,\ell} \text{R}_{i \bar{j} k \bar{\ell}} \xi^i \bar{\xi}^j \xi^k \bar{\xi}^{\ell}, \qquad \xi \in T^{1,0} X.
\end{eqnarray} These manifolds are biholomorphically isometric to the Bergman metric on $\mathbf{B}^n$ if $c<0$,  the Euclidean metric on $\mathbf{C}^n$ if $c=0$,  and the Fubini--Study metric on $\mathbf{P}^n$ if $c>0$ (see \cite{Hawley, Igusa,MogiYano}).

A long-standing conjecture predicts that a compact Hermitian manifold with constant (Chern) holomorphic sectional curvature is K\"ahler if $c \neq 0$ and Chern-flat if $c=0$.  In particular,  compact Hermitian manifolds $(X,\text{g})$ with $\text{HSC}_{\text{g}} \equiv 0$ are predicted to be compact quotients of complex Lie groups with a left-invariant metric \cite{Boothby}.  The conjecture is known for surfaces \cite{BalasGauduchon, ADM},   Chern and Strominger--Bismut K\"ahler-like metrics \cite{TangHSC, ChenZheng, RaoZheng},  complex nilmanifolds \cite{LiZheng},  and for locally conformally K\"ahler metrics with non-positive constant \cite{ChenChenNie}.  

Existing methods have all centered around forcing symmetries on the metric from derived curvature identities. While we maintain aspects of this approach here,  we also use the relationship between the holomorphic sectional curvature and the algebro-geometric properties of the canonical bundle $K_X$. In more detail, we will consider compact complex manifolds with pseudoeffective canonical bundle---i.e., we assume there is a possibly singular Hermitian metric $h$ on $K_X$ with whose curvature current $\Theta(h): = \sqrt{-1} \partial \bar{\partial} \log(h)$ is semi-positive (in the sense of currents). This forces additional structure on the metric and allows us to establish the first main theorem.

\subsection*{Theorem 1.1}\label{11}
Let $X$ be a compact complex manifold with pseudoeffective canonical bundle. Then, any Hermitian metric g on $X$ with $\text{HSC}_{\text{g}} \equiv 0$ is conformally balanced. In particular, pluriclosed metrics with $\text{HSC}_{\text{g}} \equiv 0$ are flat K\"ahler metrics and $X$ is biholomorphic to a complex torus $\mathbf{C}^n / \Lambda$. \\

\noindent We remind the reader that a Hermitian metric g is \textit{pluriclosed} if $\partial \bar{\partial} \omega_{\text{g}} =0$ and \textit{balanced} if $\text{d}\omega_{\text{g}}^{n-1} =0$, where $n : = \dim_{\mathbf{C}} X$. Note that if $X$ is Moishezon (bimeromorphic to a projective manifold) or a compact complex threefold in the Fujiki class $\mathcal{C}$ (bimeromorphic to a compact K\"ahler manifold), the curvature condition $\text{HSC}_{\text{g}} \equiv 0$ automatically implies that $K_X$ is pseudoeffective (as a consequence of the Ahlfors--Schwarz lemma and \cite{BDPP, Brunella}). We believe that the assumption on the canonical bundle is superfluous in general, and even more, predict that a compact Hermitian manifold $(X,\text{g})$ with $\text{HSC}_{\text{g}} \equiv 0$ is balanced.

For compact K\"ahler manifolds endowed with a pluriclosed metric of $\text{HSC}_{\text{g}} \equiv 0$,  we can make use of the results in \cite{BroderStanfieldSchwarz} to prove the following.

\subsection*{Theorem 1.2}\label{12}
Let $X$ be a compact K\"ahler manifold with a pluriclosed metric of $\text{HSC}_{\text{g}} \equiv 0$.  Then $(X,\text{g})$ is biholomorphically isometric to a complex torus with a flat K\"ahler metric. \\

\noindent We next consider the curvature constraint that was introduced by Yang--Zheng \cite{YangZhengRBC}.  In more detail,  motivated by the Schwarz lemma for holomorphic maps into Hermitian manifolds (see, e.g., \cite{Royden, YauSchwarzLemma, BroderSBC, BroderSurvey, BroderStanfieldSchwarz}),  Yang--Zheng introduced the \textit{real bisectional curvature} \begin{eqnarray}\label{RBCDEFN}
\text{RBC}_{\text{g}}(\zeta) & : = & \frac{1}{| \zeta |_{\text{g}}^2}   \sum_{\alpha, \beta, \gamma, \delta} \text{R}_{\alpha \bar{\beta} \gamma \bar{\delta}} \zeta^{\alpha \bar{\beta}} \zeta^{\gamma \bar{\delta}},
\end{eqnarray} where $\zeta$ is a nonnegative Hermitian $(1,1)$-tensor.  If the metric is K\"ahler(-like), the real bisectional curvature is comparable to the holomorphic sectional curvature (i.e.,  if one is signed,  the other has the same sign),  but in general,  it defines a stronger curvature.  However, it is not strong enough to control any of the Chern Ricci curvatures.  Yang--Zheng showed that compact Hermitian manifolds with $\text{RBC}_{\text{g}} \equiv 0$ are balanced with vanishing Chern Ricci curvatures.  If,  in addition,  $\dim_{\mathbf{C}} X =3$,  Zhou--Zheng \cite{ZhouZheng} showed that the metric is Chern-flat.

Using the locality technique introduced in \cite{BroderPulemotov},  we show that Hermitian metrics g with $\text{RBC}_{\text{g}} \equiv 0$ on complex manifolds in the Fujiki class $\mathcal{C}$ are flat K\"ahler metrics on a complex torus.  

\subsection*{Theorem 1.3}\label{13}
Let $X$ be a compact complex manifold in the Fujiki class $\mathcal{C}$.  Suppose that $X$ admits a Hermitian metric g with $\text{RBC}_{\text{g}} \equiv 0$.  Then $(X,\text{g})$ is biholomorphically isometric to a complex torus $\mathbf{C}^n/\Lambda$ with a flat K\"ahler metric.   \\

\noindent In what remains of the present article,  we will consider a family of `altered curvatures'.  These curvatures primarily serve the role of forcing distinguished structures on the metric if they happen to coincide with their `standard' counterparts.  The model example of an altered curvature is what we call the \textit{altered scalar curvature} $\widetilde{\text{Scal}}_{\text{g}} : = \text{g}^{i \bar{\ell}} \text{g}^{k \bar{j}} \text{R}_{i \bar{j} k \bar{\ell}}$,  whose `standard' counterpart is the Chern scalar curvature $\text{Scal}_{\text{g}} : = \text{g}^{i \bar{j}} \text{g}^{k \bar{j}} \text{R}_{i \bar{j} k \bar{\ell}}$.  The primary significance of the altered scalar curvature is expressed in an old result of Gauduchon \cite{Gauduchon},  which asserts that a compact Hermitian manifold with $\text{Scal}_{\text{g}} = \widetilde{\text{Scal}}_{\text{g}}$ is balanced.  These altered curvatures play an essential role in understanding the curvature of a Hermitian metric \cite{LiuYangRicci},  including the curvature of the Gauduchon connections  \cite{BroderStanfieldGauduchon}.  This article will introduce altered variants of the holomorphic sectional curvature, real bisectional curvature, and quadratic orthogonal bisectional curvature.

\section{Compact complex manifolds with constant holomorphic and real bisectional curvature}
\noindent In this section,  we will prove  \nameref{11},  \nameref{12},  and \nameref{13}.  Throughout this manuscript, $X$ will always denote a (connected and paracompact) complex manifold.  We exclusively denote by $n : = \dim_{\mathbf{C}} X$ the complex dimension of $X$.  If g is a Hermitian metric on $X$ and J denotes the (integrable) complex structure on $X$,  then we write $\omega_{\text{g}}(\cdot, \cdot) : = \text{g}(\text{J} \cdot, \cdot)$ for the $(1,1)$-form associated to g. The top exterior power of $\omega$ is related to the Riemannian volume form $\text{d}v_{\text{g}}$ by $\omega_{\text{g}}^n = n! \text{d}v_{\text{g}}$.

\subsection*{Reminder 2.1. }
The \textit{Chern connection} is the unique $\mathbf{C}$-linear connection on the tangent bundle $T^{1,0}X$ such that $\nabla\text{g} = \nabla \text{J} =0$ and the $(0,1)$-part coincides with the holomorphic structure $\nabla^{0,1} = \bar{\partial}$ (see,  e.g.,   \cite{GauduchonConnections, BroderStanfieldGauduchon}).   If $\text{g}_{i \bar{j}}$ denote the components of the Hermitian metric g in a local frame,  the components of the Chern torsion T are given by $\text{T}_{ij}^k = \text{g}^{k \bar{\ell}} (\partial_i \text{g}_{j \bar{\ell}} - \partial_j \text{g}_{i \bar{\ell}})$.  We write $\eta = \sum_i \eta_i dz^i = \sum_{i,k} \text{T}_{ik}^k dz^i$ for the \textit{torsion $(1,0)$-form}.

Let R denote the curvature of the Chern connection with local components $\text{R}_{i \bar{j} k \bar{\ell}}$.  There are four distinct traces of the Chern curvature tensor and, thus, four distinct Chern Ricci curvatures.  The \textit{first} and \textit{second Chern Ricci curvatures} $\text{Ric}_{i \bar{j}}^{(1)}  : =  \text{g}^{k \bar{\ell}} \text{R}_{i \bar{j} k \bar{\ell}}$ and  $\text{Ric}_{k \bar{\ell}}^{(2)}  : =  \text{g}^{i \bar{j}} \text{R}_{i \bar{j} k \bar{\ell}}$ trace to the same \textit{Chern scalar curvature} $\text{Scal}_{\text{g}}$.  The remaining \textit{third} and \textit{fourth Chern Ricci curvatures} $\text{Ric}_{i \bar{\ell}}^{(3)}  : =  \text{g}^{k \bar{j}} \text{R}_{i \bar{j} k \bar{\ell}}$ and $\text{Ric}_{k \bar{j}}^{(4)} : =  \text{g}^{i \bar{\ell}} \text{R}_{i \bar{j} k \bar{\ell}}$ trace to what we refer to as the \textit{altered Chern scalar curvature} $\widetilde{\text{Scal}}_{\text{g}}$, i.e.,   \begin{eqnarray*}
\text{Scal}_{\text{g}} \ : = \ \text{g}^{i \bar{j}} \text{Ric}_{i \bar{j}}^{(1)} \ = \ \text{g}^{k \bar{\ell}} \text{Ric}_{k \bar{\ell}}^{(2)}, \qquad \widetilde{\text{Scal}}_{\text{g}} \ : = \ \text{g}^{i \bar{\ell}} \text{Ric}_{i \bar{\ell}}^{(3)} \ = \ \text{g}^{k \bar{j}} \text{Ric}_{k \bar{j}}^{(4)}.
\end{eqnarray*} We invite the reader to consult \cite{LiuYangRicci,BroderStanfieldGauduchon} for further details.

\subsection*{Proof of \nameref{11}}
Let $X$ be a compact complex manifold with pseudoeffective canonical bundle. Let g be a Hermitian metric on $X$ with $\text{HSC}_{\text{g}} \equiv 0$. For any $\xi \in T^{1,0} X$,  we have \begin{eqnarray}\label{HSCScalars}
\int_{\mathbf{P}^{n-1}} \frac{1}{|\xi |_{\text{g}}^4} \text{R}_{i \bar{j} k \bar{\ell}} \xi^i \bar{\xi}^j \xi^k \bar{\xi}^{\ell} \omega_{\text{FS}}^{n-1} &=& \frac{1}{n(n+1)} \left( \text{Scal}_{\text{g}} + \widetilde{\text{Scal}}_{\text{g}} \right),
\end{eqnarray} where $\omega_{\text{FS}}$ denotes the Fubini--Study metric of unit volume on $\mathbf{P}^{n-1}$ (see, e.g., \cite[Lemma 4.1]{LiuSunYang});  hence,  $\text{Scal}_{\text{g}} + \widetilde{\text{Scal}}_{\text{g}} \equiv 0$.  Let $\text{g}_{\text{G}} = f_0^{\frac{1}{n-1}} \text{g}$ be the Gauduchon metric (i.e., g$_{\text{G}}$ satisfies $\partial \bar{\partial} \omega_{\text{G}}^{n-1}=0$) in the conformal class of g, where $f_0 \in \mathcal{C}^{\infty}(X,\mathbf{R})$ is a strictly positive function (such a metric always exists \cite{Gauduchon}). Let $\omega_{\text{G}}(\cdot, \cdot) : = \text{g}_{\text{G}}(\text{J} \cdot, \cdot)$. Then \begin{eqnarray*}
\int_X \text{Scal}_{\text{g}_{\text{G}}} \omega_{\text{G}}^n &=& n \int_X \text{Ric}_{\text{g}_{\text{G}}}^{(1)} \wedge \omega_{\text{G}}^{n-1} \\
&=& n \int_X \text{Ric}_{\text{g}}^{(1)} \wedge \omega_{\text{G}}^{n-1} \ = \  n \int_X f_0 \text{Ric}_{\text{g}}^{(1)} \wedge \omega_{\text{g}}^{n-1} \ = \  \int_X f_0 \text{Scal}_{\text{g}} \omega_{\text{g}}^n.
\end{eqnarray*} A similar computation (see, e.g., \cite{YangMRL}) then gives $\int_X \widetilde{\text{Scal}}_{\text{g}_{\text{G}}} \omega_{\text{G}}^n = \int_X f_0 \widetilde{\text{Scal}}_{\text{g}} \omega_{\text{g}}^n$.  Let $\eta_{\text{G}}$ denote the torsion $(1,0)$-form of $\text{g}_{\text{G}}$. We then compute  \begin{eqnarray*}
\int_X \text{Scal}_{\text{g}_{\text{G}}} \omega_{\text{G}}^n &=& \frac{1}{2} \int_X \left( \text{Scal}_{\text{g}_{\text{G}}} + \widetilde{\text{Scal}}_{\text{g}_{\text{G}}} \right) \omega_{\text{G}}^n + \frac{1}{2} \int_X \left( \text{Scal}_{\text{g}_{\text{G}}} - \widetilde{\text{Scal}}_{\text{g}_{\text{G}}} \right) \omega_{\text{G}}^n  \\
&=& \frac{1}{2} \int_X f_0 \left( \text{Scal}_{\text{g}} + \widetilde{\text{Scal}}_{\text{g}} \right) \omega_g^n + \frac{1}{2} \int_X | \eta_{\text{G}} |^2_{\text{g}_{\text{G}}} \omega_{\text{G}}^n  \\
&=&  \frac{1}{2} \int_X | \eta_{\text{G}} |^2_{\text{g}_{\text{G}}} \omega_{\text{G}}^n \ \geq \ 0.
\end{eqnarray*} Since $K_X$ is pseudoeffective, \cite[Theorem 1.3]{YangScalarCurvature} implies that $\int_X | \eta_{\text{G}} |_{\text{g}_{\text{G}}}^2 \omega_{\text{G}}^n =0$. Therefore,  $\text{g}_{\text{G}}$ is balanced, and g is conformally related to a balanced metric. This proves the first statement of \nameref{11}.  If, in addition,  g is pluriclosed,  the metric must be K\"ahler (see, e.g.,  \cite{LiuYangRicci}, Corollary 5.4). In particular,  $(X,\text{g})$ is a compact K\"ahler manifold with $\text{HSC}_{\text{g}} \equiv 0$ and thus biholomorphically isometric to a complex torus with a flat K\"ahler metric.  \hfill $\Box$ \\

\noindent An essential part of the proof of \nameref{12} is the following result.

\subsection*{Proposition 2.2}\label{24}
Let $(X,\text{g})$ be a compact Hermitian manifold with constant holomorphic sectional curvature $\text{HSC}_{\text{g}} \equiv c$, for some $c \in \mathbf{R}$.  Then \begin{eqnarray*}
\int_X \text{Scal}_{\text{g}} \omega_{\text{g}}^n &=& \frac{cn(n+1)}{2} \int_X \omega_{\text{g}}^n + \frac{1}{2} \int_X | \eta |_{\text{g}}^2 \omega_{\text{g}}^n.
\end{eqnarray*} In particular,  a pluriclosed metric with  $\text{HSC}_{\text{g}} \equiv 0$ and $\int_X \text{Scal}_{\text{g}} \omega_{\text{g}}^n =0$ is K\"ahler and the universal cover of $X$ is $\mathbf{C}^n$.  \begin{proof}
Suppose g is a Hermitian metric with $\text{HSC}_{\text{g}} \equiv c$.  Let $\{ e_i \}$ be a local unitary frame with dual coframe $\{ e^k \}$ of $(1,0)$-forms.  From the first Bianchi identity $\text{T}_{ij, \bar{\ell}}^k = \text{R}_{j\bar{\ell}i \bar{k}} - \text{R}_{i \bar{\ell} j \bar{k}}$,  we have \begin{eqnarray*}
\eta_{j, \bar{\ell}} &=& \sum_k \left( \text{R}_{j \bar{\ell} k \bar{k}} - \text{R}_{k \bar{\ell} j \bar{k}} \right), 
\end{eqnarray*} where the index after the comma indicates covariant differentiation with respect to the Chern connection.  Hence,   \begin{eqnarray*}
\sum_i \eta_{i, \bar{i}} \ = \  \sum_{i,k} \left( \text{R}_{i \bar{i} k \bar{k}} - \text{R}_{k \bar{i} i \bar{k}} \right) &=& \sum_{i \neq k} \left( \text{R}_{i \bar{i} k \bar{k}} - \text{R}_{k \bar{i} i \bar{k}} \right) \\
&=& \sum_{i<k} \left( \text{R}_{i \bar{i} k \bar{k}} + \text{R}_{k \bar{k} i \bar{i}} \right) - \sum_{i<k} \left( \text{R}_{k \bar{i} i \bar{k}} + \text{R}_{i \bar{k} k \bar{i}} \right).
\end{eqnarray*} The Balas lemma \cite{Balas} asserts that  \begin{eqnarray*}
\text{R}_{i \bar{i} k \bar{k}} + \text{R}_{k \bar{i} i \bar{k}} + \text{R}_{i \bar{k} k \bar{i}} + \text{R}_{k \bar{k} i \bar{i}} &=& 2c,
\end{eqnarray*} and therefore, \begin{eqnarray*}
\sum_{i<k} \left( \text{R}_{k \bar{i} i \bar{k}} + \text{R}_{i \bar{k} k \bar{i}} \right) &=& \sum_{i<k} \left[ 2c - (\text{R}_{i \bar{i} k \bar{k}} + \text{R}_{k \bar{k} i \bar{i}} ) \right].
\end{eqnarray*} Inserting this into the previous computation,   \begin{eqnarray*}
\sum_i \eta_{i, \bar{i}} &=& 2 \sum_{i<k} \left( \text{R}_{i \bar{i} k \bar{k}} + \text{R}_{k \bar{k} i \bar{i}} \right) - cn(n-1) \\
&=& 2 \left( \text{Scal}_{\text{g}} - nc \right) - cn(n-1) \ = \ 2 \text{Scal}_{\text{g}} - cn(n+1).
\end{eqnarray*} The torsion $(1,0)$-form satisfies $
\partial \omega_{\text{g}}^{n-1} = - \eta \wedge \omega_{\text{g}}^{n-1}$.  Hence,  \begin{eqnarray*}
\bar{\partial} \partial \omega_{\text{g}}^{n-1} &=& - \bar{\partial} \eta \wedge  \omega_{\text{g}}^{n-1} + \eta \wedge \bar{\partial} \omega_{\text{g}}^{n-1} \ = \ - \bar{\partial} \eta \wedge \omega_{\text{g}}^{n-1} - \eta \wedge \bar{\eta} \wedge \omega_{\text{g}}^{n-1}.
\end{eqnarray*} Hence,  by integrating,   \begin{eqnarray*}
\int_X  | \eta |_{\text{g}}^2 \omega_{\text{g}}^n \ = \  \int_X \left( \sum_i \eta_{i,\bar{i}} \right) \omega_{\text{g}}^n &=& 2 \int_X \text{Scal}_{\text{g}} \omega_{\text{g}}^n - cn(n+1) \int_X \omega_{\text{g}}^n.
\end{eqnarray*} 
\end{proof}

\subsection*{Proof of \nameref{12}}
Let $X$ be a compact K\"ahler manifold with a pluriclosed metric of $\text{HSC}_{\text{g}} \equiv 0$. Applying \cite[Theorem 3.2]{BroderStanfieldSchwarz}, we observe that $K_X$ is nef (i.e., for any $\varepsilon>0$, there exists a smooth Hermitian metric $h$ with curvature form $\Theta(h) : = \sqrt{-1} \partial \bar{\partial} \log(h) \geq - \varepsilon \omega_0$, where $\omega_0$ is a fixed reference form). In particular, $K_X$ is pseudoeffective. As in the proof of \nameref{11}, the vanishing of the holomorphic sectional curvature implies $\text{Scal}_{\text{g}} + \widetilde{\text{Scal}}_{\text{g}} \equiv 0$. Since g is pluriclosed, we also have the following identity on the curvature and torsion of g: \begin{eqnarray}\label{PluriclosedIdentity}
\text{R}_{i \bar{j} k \bar{\ell}} + \text{R}_{k \bar{\ell} i \bar{j}} - \text{R}_{i \bar{\ell} k \bar{j}} - \text{R}_{k \bar{j} i \bar{\ell}} &=& \text{T}_{ik}^p \overline{\text{T}_{j\ell}^q} \text{g}_{p \bar{q}}.
\end{eqnarray} Using the metric to trace over $i,j$ and $k,\ell$, we see that pluriclosed metrics satisfy $\text{Scal}_{\text{g}} - \widetilde{\text{Scal}}_{\text{g}} = \frac{1}{2} | \text{T} |^2 \geq 0$. In particular, pluriclosed metrics with $\text{HSC}_{\text{g}} \equiv 0$ satisfy $\text{Scal}_{\text{g}} \geq 0$.  Now proceed by contradiction and suppose that g is not balanced.  \nameref{24} with $c=0$ implies that g has positive total Chern scalar curvature. Let $ \text{g}_{\text{G}} = f_0^{\frac{1}{n-1}} \text{g} $ denote the Gauduchon metric in the conformal class of g. From the proof of \nameref{11}, we see that \begin{eqnarray*}
\int_X \text{Scal}_{\text{g}_{\text{G}}} \omega_{\text{G}}^n &=& \int_X f_0 \text{Scal}_{\text{g}} \omega_{\text{g}}^n \ \geq \ \left( \min_{x \in X} f_0(x) \right) \int_X \text{Scal}_{\text{g}} \omega_{\text{g}}^n \ > \ 0.
\end{eqnarray*} Hence, $X$ admits a Gauduchon metric with positive total Chern scalar curvature. From \cite[Theorem 1.3]{YangScalarCurvature},  this implies that $K_X$ is not pseudoeffective and thus violates the fact that $K_X$ is nef. It follows that g is balanced. A Hermitian metric that is both pluriclosed and balanced is K\"ahler. Therefore,  $(X,\text{g})$ is a compact K\"ahler manifold with $\text{HSC}_{\text{g}} \equiv 0$; such manifolds are completely classified by \cite{Hawley, Igusa, MogiYano}.\hfill $\Box$

\hfill

\noindent Recall that the Chern scalar curvature is said to be \textit{quasi-positive} if $\text{Scal}_{\text{g}} \geq 0$ on $X$ and $\text{Scal}_{\text{g}}(x_0) > 0$ for some $x_0 \in X$. From the proof of \nameref{12}, we also observe the following extension of a result in \cite{YangScalarCurvature}.

\subsection*{Corollary 2.3}
Let $X$ be a compact complex manifold. The following are equivalent: \begin{itemize}
    \item[(i)] $X$ admits a Hermitian metric with quasi-positive Chern scalar curvature.
    \item[(ii)] $X$ admits a Gauduchon metric with positive total Chern scalar curvature. 
\end{itemize} In particular, a compact Hermitian manifold with quasi-positive Chern scalar curvature has negative Kodaira dimension. \begin{proof}
    The equivalence of (i) and (ii) follows from the proof of \nameref{12}, as we already mentioned. Hence, if $X$ has a Hermitian metric with quasi-positive Chern scalar curvature, we can find a Gauduchon metric with positive total Chern scalar curvature. From \cite[Lemma 3.2]{YangMRL}, there is a Hermitian metric with positive Chern scalar curvature. Proceed by contradiction and suppose that the Kodaira dimension is nonnegative. Then there exists a holomorphic section $\sigma \in H^0(X, K_X^{\otimes \ell})$ for some $\ell \in \mathbf{N}$. Hence, a simple calculation gives \begin{eqnarray*}
        \Delta_{\text{g}} | \sigma |^2_{\text{g}} &=& | \nabla \sigma |^2 + \ell \text{Scal}_{\text{g}} | \sigma |_{\text{g}}^2 \ > \ | \nabla \sigma |^2.
    \end{eqnarray*} The divergence theorem then yields the desired contradiction.
\end{proof}

\subsection*{Remark 2.4}\label{25}
Let us remark that an essential feature of the proof of \nameref{12} is that by using the Schwarz lemma \cite{BroderStanfieldSchwarz} to prove that the canonical bundle $K_X$ is nef,  we can extract more precise information out of Yang's result on the image of the total Chern scalar curvature functional. The K\"ahler assumption is only used in order to apply the Wu--Yau theorem given in \cite{BroderStanfieldSchwarz}. While it is expected that a compact Hermitian manifold $(X,\text{g})$ with $\text{HSC}_{\text{g}} \leq 0$ has nef canonical bundle,  the existing strategies \cite{WuYau1, TosattiYang, YangZhengRBC, BroderStanfieldSchwarz},  this presently remains out of reach.  \\

\noindent Let us mention that \nameref{12} extends to complex manifolds $X$ is in the Fujiki class $\mathcal{C}$,  i.e.,  $X$ is bimeromorphic to a compact K\"ahler manifold,  endowed with a pluriclosed metric of constant $\text{HSC}_{\text{g}} \equiv 0$. Indeed,  since $X$ has no rational curves,  i.e., holomorphic maps $\mathbf{P}^1 \to X$,  we know that $X$ is K\"ahler (see, e.g.,  \cite{BiswasMcKay}).  

If we remove the pluriclosed assumption and consider Hermitian metrics with constant real bisectional curvature (defined in \eqref{RBCDEFN}),  then we can make use of the locality technique introduced in \cite{BroderPulemotov} to prove \nameref{13}.

\subsection*{Proof of \nameref{13}}
We first note that the Schwarz lemma \cite{YangZhengRBC, BroderStanfieldSchwarz,BroderSurvey} shows that if $X$ has a Hermitian metric with $\text{RBC}_{\text{g}} \equiv 0$,  then every holomorphic $\mathbf{P}^1 \to X$ is constant.  Hence,  from \cite{BiswasMcKay},  $X$ must be K\"ahler.  Yang--Zheng \cite{YangZhengRBC} showed that a Hermitian metric with $\text{RBC}_{\text{g}} \equiv 0$ is balanced with vanishing Chern Ricci curvatures.  In particular, the canonical bundle $K_X$ is holomorphically torsion,  i.e.,  $K_X^{\otimes \ell} \simeq \mathcal{O}_X$ for some $\ell \in \mathbf{N}$.  From Yau's resolution of the Calabi conjecture \cite{Yau1976},  $X$ admits a Ricci-flat K\"ahler metric $\hat{\text{g}}$.  Using the technique introduced in \cite{BroderPulemotov}, we apply the Schwarz lemma to show that the Chern connections of g and $\hat{\text{g}}$ coincide.  In more detail,  let $\text{id} : (X,\hat{\text{g}}) \to (X, \text{g})$ denote the identity map.  Observe that $(\nabla \partial \text{id})_{ij}^k = \hat{\Gamma}_{ij}^k - \Gamma_{ij}^k$,  where $\hat{\Gamma}_{ij}^k = \hat{\text{g}}^{k \bar{\ell}} \partial_i \hat{\text{g}}_{j \bar{\ell}}$ and $\Gamma_{ij}^k = \text{g}^{k \bar{\ell}} \partial_i \text{g}_{j \bar{\ell}}$ denote the local coefficients for the Chern connection of $\hat{\text{g}}$ and g,  respectively.  Hence,  to show that $\hat{\text{g}}$ and g have the same Chern connection, it suffices to show that $| \nabla \partial \text{id} |^2 =0$.  To this end,  we compute \begin{eqnarray*}
\Delta_{\hat{\text{g}}} \text{tr}_{\hat{\text{g}}}(\text{g}) & = & | \nabla \partial \text{id} |^2 + \text{Ric}(\hat{\text{g}})_{k \bar{\ell}} \hat{\text{g}}^{p \bar{\ell}} \hat{\text{g}}^{k \bar{q}} \text{g}_{p \bar{q}} - \text{R}_{i \bar{j} k \bar{\ell}} \hat{\text{g}}^{i \bar{j}} \hat{\text{g}}^{k \bar{\ell}},
\end{eqnarray*} where $\text{Ric}(\hat{\text{g}})$ is the  Ricci curvature of $\hat{\text{g}}$ and $\text{R}$ is Chern curvature tensor of g.  Since $\hat{\text{g}}$ is Ricci-flat K\"ahler and the real bisectional curvature of g vanishes,  the divergence theorem implies \begin{eqnarray*}
0 \ = \ \int_X \Delta_{\hat{\text{g}}} \text{tr}_{\hat{\text{g}}}(\text{g}) \omega_{\text{g}}^n & =  & \int_X | \nabla \partial \text{id} |^2 \omega_{\text{g}}^n.
\end{eqnarray*} Hence,  $| \nabla \partial \text{id} |^2=0$ and the Chern connections of g and $\hat{\text{g}}$ coincide.  The K\"ahler condition is expressed by the vanishing of the torsion of the Chern connection.  Since the Hermitian metric g has the same Chern connection as the K\"ahler metric $\hat{\text{g}}$,  it follows that g is K\"ahler.  A K\"ahler metric with $\text{RBC}_{\text{g}} \equiv 0$ has $\text{HSC}_{\text{g}} \equiv 0$ and hence,  $(X, \text{g})$ is (biholomorphically isometric to) a complex torus with a flat K\"ahler metric. \hfill $\Box$ 

\hfill

\noindent A further corollary of \nameref{13} and the locality technique of \cite{BroderPulemotov} is the following rigidity statement for second Chern Ricci-flat metrics.

\subsection*{Corollary 2.5}
Let $X$ be a compact complex manifold in the Fujiki class $\mathcal{C}$.  Suppose $X$ admits a pluriclosed metric with $\text{RBC}_{\text{g}} \equiv 0$ and a Hermitian metric $\hat{\text{g}}$ with $\text{Ric}_{\hat{\text{g}}}^{(2)}=0$.  Then$X$ is a complex torus, and both g and $\hat{\text{g}}$ are flat K\"ahler metrics.  

\section{The Altered Curvatures}
\noindent In this section,  we describe the altered curvatures in greater detail. Recall that an altered curvature is a curvature tensor that forces distinguished structure on a metric when compared to its `standard' counterpart. The model examples of altered curvatures are the third and fourth Chern Ricci curvatures and their trace,  which we refer to as the altered scalar curvature. We first consider an altered variant of the real bisectional curvature.

\subsection*{Definition 3.1}
Let $(X,\text{g})$ be a Hermitian manifold. For a nonnegative Hermitian $(1,1)$-tensor $\zeta$,  the \textit{altered real bisectional curvature} is defined \begin{eqnarray}
\widetilde{\text{RBC}}_{\text{g}}(\zeta) &:=& \frac{1}{| \zeta |_{\text{g}}^2}  \sum_{\alpha, \beta, \gamma, \delta} \text{R}_{\alpha \bar{\beta} \gamma \bar{\delta}} \zeta^{\alpha \bar{\delta}} \zeta^{\gamma \bar{\beta}}.
\end{eqnarray} 

\noindent The altered real bisectional curvature controls the holomorphic sectional curvature and the altered Chern scalar curvature. Indeed,  for any unit $(1,0)$-tangent vector $\xi \in T^{1,0}X$, we may choose a unitary frame $\{ e_1, ..., e_n \}$ such that $\xi = \lambda e_1$ for some $\lambda \in \mathbf{C}$.  Hence,  if we choose $\zeta^{\alpha \bar{\beta}} = \delta^{1 \alpha} \delta^{1 \beta}$,  we see that \begin{eqnarray*}
\widetilde{\text{RBC}}_{\text{g}}(\zeta) &=& \sum_{\alpha, \beta, \gamma, \delta} \text{R}_{\alpha \bar{\delta} \gamma \bar{\beta}} \delta^{\alpha \bar{\beta}} \delta^{\gamma \bar{\delta}} \ = \ \text{R}_{1 \bar{1} 1 \bar{1}} \ = \ \text{HSC}_{\text{g}}(\xi).
\end{eqnarray*} Further,  choosing $\zeta^{\alpha \bar{\beta}} = \frac{1}{\sqrt{n}} \delta_{\alpha \beta}$ we have $\widetilde{\text{RBC}}_{\text{g}}(\zeta) = \widetilde{\text{Scal}}_{\text{g}}$.  

The following result provides an analog of what is known for Hermitian metrics with constant real bisectional curvature (but curiously opposite in sign).

\subsection*{Theorem 3.2}\label{31}
Let $(X,\text{g})$ be a compact Hermitian manifold with $\widetilde{\text{RBC}}_{\text{g}} \equiv c$,  for some $c \in \mathbf{R}$.  Then $c \geq 0$,  and if $c=0$,  the metric g is balanced with vanishing Chern Ricci curvatures.   \begin{proof}

Suppose g is a Hermitian metric with $\widetilde{\text{RBC}}_{\text{g}} \equiv c$.  For any local frame and nonnegative Hermitian $(1,1)$-tensor $\zeta$,  we have \begin{eqnarray*}
\text{R}_{i \bar{j} k \bar{\ell}} \zeta^{i \bar{\ell}} \zeta^{k \bar{j}} &=& c \text{tr}(\zeta^2). 
\end{eqnarray*} In particular, \begin{eqnarray}\label{RBCConstantFormulae}
\text{R}_{i \bar{i} k \bar{k}} + \text{R}_{k \bar{k} i \bar{i}} \ = \ 2c, \qquad \text{R}_{i \bar{j} k \bar{\ell}} + \text{R}_{k \bar{\ell} i \bar{j}} \ = \ 0.
\end{eqnarray} Hence, in a similar manner to the proof of \nameref{24}, we compute \begin{eqnarray*}
\sum_i \eta_{i, \bar{i}} \ = \  \sum_{i,k} \left( \text{R}_{i \bar{i} k \bar{k}} - \text{R}_{k \bar{i} i \bar{k}} \right) &=& \sum_{i \neq k} \left( \text{R}_{i \bar{i} k \bar{k}} - \text{R}_{k \bar{i} i \bar{k}} \right) \\
&=& \sum_{i \neq k} \left( 2c - \text{R}_{k \bar{k} i \bar{i}} + \text{R}_{i \bar{k} k \bar{i}} \right) \\
&=& 2c n(n-1)- \sum_{i \neq k} (\text{R}_{k \bar{k} i \bar{i}} - \text{R}_{i \bar{k} k \bar{i}} ) \\
&=& 2c n (n-1) - \sum_i \eta_{i, \bar{i}}.
\end{eqnarray*} It follows that \begin{eqnarray*}
\int_X | \eta |_{\text{g}}^2 \omega_{\text{g}}^n &=& \int_X \left( \sum_i \eta_{i,\bar{i}} \right) \omega_{\text{g}}^n \ = \ cn (n-1),
\end{eqnarray*} and in particular,  $c \geq 0$.  If $c=0$,  then the metric g is balanced and thus $\partial^{\ast}_{\text{g}} \omega_{\text{g}} =0$ and $\bar{\partial}_{\text{g}}^{\ast} \omega_{\text{g}} =0$.  From \eqref{RBCConstantFormulae},  we see that the first Chern Ricci curvature $\text{Ric}_{i \bar{j}}^{(1)} = \text{g}^{k \bar{\ell}}\text{R}_{i \bar{j} k \bar{\ell}}$ vanishes if $c=0$.  Further,  since $\text{Ric}_{\omega_{\text{g}}}^{(3)} = \text{Ric}_{\omega_{\text{g}}}^{(1)} - \partial \partial^{\ast}_{\text{g}} \omega_{\text{g}}$ and $\text{Ric}_{\omega_{\text{g}}}^{(4)} = \text{Ric}_{\omega_{\text{g}}}^{(1)} - \bar{\partial} \bar{\partial}^{\ast}_{\text{g}} \omega_{\text{g}}$ (see, e.g., \cite{LiuYangRicci}),  this implies the vanishing of the third and fourth Chern Ricci curvatures. From the second equation in \eqref{RBCConstantFormulae},  this also implies the vanishing of the second Chern Ricci curvature.  
\end{proof}

\noindent Consideration of the altered real bisectional curvature is motivated by the fact that the holomorphic sectional curvature is comparable to the sum of the real bisectional curvature and altered real bisectional curvature.  

\subsection*{Definition 3.3}
For nonnegative Hermitian $(1,1)$-tensors $\zeta$,  we define the \textit{altered holomorphic sectional curvature} \begin{eqnarray*}
\widetilde{\text{HSC}}_{\text{g}}(\zeta) & := & \frac{1}{| \zeta |_{\text{g}}^2} \sum_{\alpha, \beta, \gamma, \delta} \left( \text{R}_{\alpha \bar{\beta} \gamma \bar{\delta}} + \text{R}_{\alpha \bar{\delta} \gamma \bar{\beta}} \right) \zeta^{\alpha \bar{\beta}} \zeta^{\gamma \bar{\delta}}.
\end{eqnarray*} This curvature function is considered implicitly in \cite{YangZhengRBC} and is comparable to the familiar holomorphic sectional curvature in the sense that if one is signed,  the other admits the same sign (see,  e.g., \cite[p. 5]{YangZhengRBC}). These two holomorphic sectional curvatures are not scalar multiples of each other, however.

In fact,  the following result shows that constant holomorphic sectional curvature does not imply constant altered holomorphic sectional curvature, in general.

\subsection*{Proposition 3.4}\label{32}
Let $(X, \text{g})$ be a Hermitian manifold with constant holomorphic sectional curvature $\text{HSC}_{\text{g}} \equiv c$, for some $c \in \mathbf{R}$.  Then the altered holomorphic sectional curvature satisfies \begin{eqnarray*}
\widetilde{\text{HSC}}_{\text{g}}(\zeta) &=& \frac{c}{| \zeta |^2} \left( 1 + \sum_{i,k=1}^n \zeta^{i \bar{k}} \zeta^{k \bar{i}} \right).
\end{eqnarray*} In particular,  $\text{HSC}_{\text{g}} \equiv c$ and $ \widetilde{\text{HSC}}_{\text{g}} \equiv c$ if and only if $c=0$.  \begin{proof}
Suppose $\text{HSC}_{\text{g}} \equiv c$ for some $c \in \mathbf{R}$.  In any local unitary frame,  the Balas lemma \cite{Balas} yields \begin{eqnarray*}
\text{R}_{i \bar{i} k \bar{k}} + \text{R}_{k \bar{i} i \bar{k}} + \text{R}_{i \bar{k} k \bar{i}} + \text{R}_{k \bar{k} i \bar{i}} &=& 2c.
\end{eqnarray*} Let $\zeta$ be a nonnegative Hermitian $(1,1)$-tensor.  Then \begin{eqnarray*}
\sum_{i\neq k} \left( \text{R}_{i \bar{i} k \bar{k}} + \text{R}_{k \bar{i} i \bar{k}} + \text{R}_{i \bar{k} k \bar{i}} + \text{R}_{k \bar{k} i \bar{i}} \right) \zeta^{i \bar{k}} \zeta^{k \bar{i}} \ = \ 2 \sum_{i \neq k} \left( \text{R}_{i \bar{i} k \bar{k}} + \text{R}_{i \bar{k} k \bar{i}} \right) \zeta^{i \bar{k}} \zeta^{k \bar{i}} &=& 2c \sum_{i \neq k} \zeta^{i \bar{k}} \zeta^{k \bar{i}}.
\end{eqnarray*}  Hence,  from $\text{R}_{i \bar{i} i \bar{i}} =c$, we have  \begin{eqnarray*}
\sum_{i,k=1}^n \left( \text{R}_{i \bar{i} k \bar{k}} + \text{R}_{i \bar{k}k \bar{i}} \right) \zeta^{i \bar{k}} \zeta^{k \bar{i}} &=& c \sum_{i \neq k} \zeta^{i \bar{k}} \zeta^{k \bar{i}} + 2c \sum_k (\zeta^{k \bar{k}})^2 \ = \ \frac{c}{| \zeta |^2} \left( 1 + \sum_{i,k=1}^n \zeta^{i \bar{k}} \zeta^{k \bar{i}} \right).\end{eqnarray*} This last expression also shows that if,  in addition,  the altered holomorphic sectional curvature is constant,  the constant is necessarily zero.
\end{proof}

\subsection*{Remark 3.5}\label{33}
We confess that this is the one instance where our diligent use of the altered terminology is abandoned. More appropriate terminology would be given by employing the prefix quasi,  pseudo,  almost, or semi. These are in conflict with other strictly defined rules; however, at the present time,  we regrettably stick with this terminology.  \\

\noindent The last of the curvatures that we will consider here is the \textit{quadratic orthogonal bisectional curvature} (QOBC) \begin{eqnarray}
\text{QOBC}_{\text{g}}  :  \mathcal{F}_X \times \mathbf{R}^n_+ \to \mathbf{R},  \qquad \text{QOBC}_{\text{g}}(\xi) \ = \ \frac{1}{| \xi |_{\text{g}}^2} \sum_{\alpha, \gamma=1}^n \text{R}_{\alpha \bar{\alpha} \gamma \bar{\gamma}}(\xi_{\alpha} - \xi_{\gamma})^2, 
\end{eqnarray} where $\mathcal{F}_X$ denotes the unitary frame bundle of $X$,  was first explicitly introduced by Wu--Yau--Zheng \cite{WuYauZheng}.  The condition $\text{QOBC}_{\text{g}} \geq 0$ is then understood to mean for all $\xi \in \mathbf{R}_+^n$ and all unitary frames.   Wu--Yau--Zheng \cite[Theorem 1]{WuYauZheng} showed that a compact K\"ahler manifold with $\text{QOBC}_{\text{g}} \geq 0$ has the property that every nef class admits a semipositive representative.  Note that the QOBC first appears implicitly in the paper of Bishop--Goldberg \cite{BishopGoldberg} as the Weitzenb\"ock curvature operator acting on real $(1,1)$-forms,  and the first invariant expression of the QOBC was given by Ni--Tam \cite[(A.5)]{NiTam} (c.f., \cite{BroderQOBC}).   Chau--Tam \cite{ChauTam} and Niu \cite{Niu} showed that a compact K\"ahler manifold $(X,\text{g})$ with nonnegative QOBC has nonnegative scalar curvature.  Based on the results in \cite{BroderQOBC}, it would be interesting to give a combinatorial proof of this statement.  For a general Hermitian metric,  we have the following.

\subsection*{Proposition 3.6}\label{36}
Let $(X,\text{g})$ be a Hermitian manifold with $\text{QOBC}_{\text{g}} \geq 0$.  Then for any unitary pair of $(1,0)$-tangent vectors $\xi, \nu$,  we have \begin{eqnarray*}
\text{Ric}_{\text{g}}^{(1)}(\xi, \bar{\xi}) + \text{Ric}_{\text{g}}^{(1)}(\nu, \bar{\nu}) + \text{Ric}_{\text{g}}^{(2)}(\xi, \bar{\xi}) + \text{Ric}_{\text{g}}^{(2)}(\nu, \bar{\nu}) & \geq & 2 \left(\text{R}(\xi, \bar{\nu}, \nu, \bar{\xi}) + \text{R}(\nu, \bar{\xi}, \xi, \bar{\nu}) \right) .
\end{eqnarray*} In particular, in any local unitary frame,  the scalar curvature satisfies \begin{eqnarray*}
\text{Scal}_{\text{g}} & \geq & \frac{1}{n-1} \sum_{1 \leq k < \ell \leq n} \left( \text{R}_{k \bar{\ell} \ell \bar{k}} + \text{R}_{\ell \bar{k} k \bar{\ell}} \right).
\end{eqnarray*} \begin{proof}
Let g be a Hermitian metric with $\text{QOBC}_{\text{g}} \geq 0$.  Then for any $\xi = (\xi_1, ..., \xi_n) \in \mathbf{R}^n$ and any local unitary frame,  we have $\sum_{i,k=1}^n \text{R}_{i \bar{i} k \bar{k}}(\xi_i - \xi_k)^2 \geq 0$.  For distinct indices $j, k, \ell$, set $\xi_k =0$, $\xi_{\ell} = 2$, and $\xi_j=1$. This gives \begin{eqnarray}\label{eqn1}
&& 4 \text{R}_{k \overline{k} \ell \overline{\ell}} + 4 \text{R}_{\ell \overline{\ell} k \overline{k}} + \sum_{j \neq k, j \neq \ell} (\text{R}_{k \overline{k} j \overline{j}} + \text{R}_{j \overline{j} k \overline{k}} + \text{R}_{\ell \overline{\ell} j \overline{j}} + \text{R}_{j \overline{j} \ell \overline{\ell}})  \nonumber \\
&& \hspace*{1cm} =  4 (\text{R}_{k \overline{k} \ell \overline{\ell}} + \text{R}_{\ell \overline{\ell} k \overline{k}}) + \sum_{j \neq k, j \neq \ell} (\text{R}_{k \overline{k} j \overline{j}} + \text{R}_{j \overline{j} k \overline{k}} + \text{R}_{\ell \overline{\ell} j \overline{j}} + \text{R}_{j \overline{j} \ell \overline{\ell}}) \ \geq \ 0.
\end{eqnarray} Let $f_k = \frac{1}{\sqrt{2}}(e_k - e_{\ell})$, $f_{\ell} = \frac{1}{\sqrt{2}}(e_k + e_{\ell})$ and $f_j = e_j$. Then \eqref{eqn1} in this frame gives \begin{eqnarray*}
&& \text{R}(e_k - e_{\ell}, \overline{e_k - e_{\ell}}, e_k + e_{\ell}, \overline{e_k + e_{\ell}}) + \text{R}(e_k + e_{\ell}, \overline{e_k + e_{\ell}}, e_k - e_{\ell}, \overline{e_k - e_{\ell}})\\
&& \hspace*{1.5cm}  + \frac{1}{2} \sum_{j \neq k, j \neq \ell} ( (\text{R}(e_k - e_{\ell}, \overline{e_k - e_{\ell}}, e_j, \overline{e_j}) + \text{R}(e_j, \overline{e_j}, e_k - e_{\ell}, \overline{e_k - e_{\ell}}) ) \\
&& \hspace*{3cm}  + \frac{1}{2} \sum_{j \neq k, j \neq \ell} \text{R}(e_k + e_{\ell}, \overline{e_k + e_{\ell}}, e_j, \overline{e_j}) + \text{R}(e_j, \overline{e_j}, e_k + e_{\ell}, \overline{e_k + e_{\ell}}) \\
&=& \text{R}_{k \overline{k} k \overline{k}} + \text{R}_{\ell \overline{\ell} \ell \overline{\ell}} + \text{R}_{k \overline{k} \ell \overline{\ell}} - \text{R}_{k \overline{\ell} \ell \overline{k}} - \text{R}_{k \overline{\ell} k \overline{\ell}} + \text{R}_{\ell \overline{\ell} k \overline{k}} - \text{R}_{\ell \overline{k} \ell \overline{k}} - \text{R}_{\ell \overline{k} k \overline{\ell}} \\
&& + \text{R}_{k \overline{k} k \overline{\ell}} + \text{R}_{k \overline{k} \ell \overline{k}} - \text{R}_{k \overline{\ell} k \overline{k}} - \text{R}_{\ell \overline{k} k \overline{k}} - \text{R}_{k \overline{\ell} \ell \overline{\ell}} + \text{R}_{\ell \overline{\ell} \ell \overline{k}} + \text{R}_{\ell \overline{\ell} k \overline{\ell}} - \text{R}_{\ell \overline{k} \ell \overline{\ell}} \\
&& + \text{R}_{k \overline{k} k \overline{k}} + \text{R}_{\ell \overline{\ell} \ell \overline{\ell}} + \text{R}_{k \overline{k} \ell \overline{\ell}} - \text{R}_{k \overline{\ell} \ell \overline{k}} - \text{R}_{k \overline{\ell} k \overline{\ell}} + \text{R}_{\ell \overline{\ell} k \overline{k}} - \text{R}_{\ell \overline{k} \ell \overline{k}} - \text{R}_{\ell \overline{k} k \overline{\ell}} \\
&& - \text{R}_{k \overline{k } k \overline{\ell}} - \text{R}_{k \overline{k} \ell \overline{k}} + \text{R}_{k \overline{\ell} k \overline{k}} + \text{R}_{\ell \overline{k} k \overline{k}} + \text{R}_{k \overline{\ell} \ell \overline{\ell}} - \text{R}_{\ell \overline{\ell} \ell \overline{k}} - \text{R}_{\ell \overline{\ell} k \overline{\ell}} + \text{R}_{\ell \overline{k} \ell \overline{\ell}} \\
&& + \frac{1}{2} \sum_{j \neq k, j \neq \ell} \left( \text{R}_{k \overline{k} j \overline{j}} + \text{R}_{\ell \overline{\ell} j \overline{j}} - \text{R}_{k \overline{\ell} j \overline{j}} - \text{R}_{\ell \overline{k} j \overline{j}} + \text{R}_{j \overline{j} k \overline{k}} + \text{R}_{j \overline{j} \ell \overline{\ell}} - \text{R}_{j \overline{j} k \overline{\ell}} - \text{R}_{j \overline{j} \ell \overline{k} } \right) \\
&& + \frac{1}{2} \sum_{ j \neq k, j \neq \ell} \left( \text{R}_{k \overline{j} j \overline{j}} + \text{R}_{\ell \overline{\ell} j \overline{j}} + \text{R}_{k \overline{\ell} j \overline{j}} + \text{R}_{\ell \overline{k} j \overline{j}} + \text{R}_{j \overline{j} k \overline{k}} + \text{R}_{j \overline{j} \ell \overline{\ell}} + \text{R}_{j \overline{j} k \overline{\ell}} + \text{R}_{j \overline{j} \ell \overline{k}}   \right) \\
&=& \text{R}_{k \overline{k} k \overline{k}} + \text{R}_{\ell \overline{\ell} \ell \overline{\ell}} + \text{R}_{k \overline{k} \ell \overline{\ell}}+ \text{R}_{\ell \overline{\ell} k \overline{k}} - \text{R}_{k \overline{\ell} \ell \overline{k}} - \text{R}_{k \overline{\ell} k \overline{\ell}} - \text{R}_{\ell \overline{k} \ell \overline{k}} - \text{R}_{\ell \overline{k} k \overline{\ell}} \\
&& \hspace*{5cm} + \frac{1}{2} \sum_{j \neq k, j \neq \ell} \left( \text{R}_{k \overline{k} j \overline{j}} + \text{R}_{\ell \overline{\ell} j \overline{j}} + \text{R}_{j \overline{j} k \overline{k}} + \text{R}_{j \overline{j} \ell \overline{\ell}}  \right) \ \geq \ 0.
\end{eqnarray*} Similarly, setting $f_k = \frac{1}{\sqrt{2}}(e_k - \sqrt{-1} e_{\ell})$, $f_{\ell} = \frac{1}{\sqrt{2}}(e_k + \sqrt{-1} e_{\ell})$ and $f_j = e_j$ gives \begin{eqnarray*}
&& \text{R}_{k \overline{k} k \overline{k}} + \text{R}_{\ell \overline{\ell} \ell \overline{\ell}} + \text{R}_{k \overline{k} \ell \overline{\ell}} + \text{R}_{\ell \overline{\ell} k \overline{k}} - \text{R}_{k \overline{\ell} \ell\overline{k}} + \text{R}_{k \overline{\ell} k \overline{\ell}} + \text{R}_{\ell \overline{k} \ell \overline{k}} - \text{R}_{\ell \overline{k} k \overline{\ell}} \\
&& \hspace*{4cm} + \frac{1}{2} \sum_{j \neq k, j \neq \ell} \left( \text{R}_{k \overline{k} j \overline{j}} + \text{R}_{\ell \overline{\ell} j \overline{j}} + \text{R}_{j \overline{j} k \overline{k}} + \text{R}_{j \overline{j} \ell \overline{\ell}}    \right)  \ \geq \ 0.
\end{eqnarray*} Adding these equations together, we get \begin{eqnarray*}
\text{R}_{k \overline{k} k \overline{k}} + \text{R}_{\ell \overline{\ell} \ell \overline{\ell}} + \text{R}_{k \overline{k} \ell \overline{\ell}} + \text{R}_{\ell \overline{\ell} k \overline{k}} - \text{R}_{k \overline{\ell} \ell \overline{k}} - \text{R}_{\ell \overline{k} k \overline{\ell}} + \sum_{j \neq k, j \neq \ell} \left( \text{R}_{k \overline{k} j \overline{j}} + \text{R}_{\ell \overline{\ell} j \overline{j}} + \text{R}_{j \overline{j} k \overline{k}} + \text{R}_{j \overline{j} \ell \overline{\ell}} \right) & \geq & 0.
\end{eqnarray*} Observe that \begin{eqnarray*}
&& \text{Ric}_{k \overline{k}}^{(1)} + \text{Ric}_{\ell \overline{\ell}}^{(1)} + \text{Ric}_{k \overline{k}}^{(2)} + \text{Ric}_{\ell \overline{\ell}}^{(2)} \\
&=& \text{R}_{k \overline{k} k \overline{k}} + \text{R}_{\ell \overline{\ell} \ell \overline{\ell}} + \text{R}_{k \overline{k} \ell \overline{\ell}} + \text{R}_{\ell \overline{\ell} k \overline{k}} + \sum_{j \neq k, j \neq \ell} (\text{R}_{k \overline{k} j \overline{j}} + \text{R}_{\ell \overline{\ell} j \overline{j}}) \\
&& \hspace*{1cm} + \text{R}_{k \overline{k} k \overline{k}} + \text{R}_{\ell \overline{\ell} \ell \overline{\ell}} + \text{R}_{\ell \overline{\ell} k \overline{k}} + \text{R}_{k \overline{k} \ell \overline{\ell}} + \sum_{j \neq k, j \neq \ell} ( \text{R}_{j \overline{j} k \overline{k}} + \text{R}_{j \overline{j} \ell \overline{\ell}}) \\
&=& 2 \left(\text{R}_{k \overline{k} k \overline{k}} + \text{R}_{\ell \overline{\ell} \ell \overline{\ell}} + \text{R}_{k \overline{k} \ell \overline{\ell}} + \text{R}_{\ell \overline{\ell} k \overline{k}} \right) + \sum_{j \neq k, j \neq \ell} (\text{R}_{k \overline{k} j \overline{j}} + \text{R}_{\ell \overline{\ell} j \overline{j}} + \text{R}_{j \overline{j} k \overline{k}} + \text{R}_{j \overline{j} \ell \overline{\ell}}).
\end{eqnarray*} Hence, for $k \neq \ell$, we have \begin{eqnarray*}
\text{Ric}_{k \overline{k}}^{(1)} + \text{Ric}_{\ell \overline{\ell}}^{(1)} + \text{Ric}_{k \overline{k}}^{(2)} + \text{Ric}_{\ell \overline{\ell}}^{(2)} & \geq & 2(\text{R}_{k \overline{\ell} \ell \overline{k}} + \text{R}_{\ell \overline{k} k \overline{\ell}}).
\end{eqnarray*} For the statement concerning the scalar curvature, we observe that \begin{eqnarray*}
2 \text{Scal}_{\text{g}} &=& \frac{1}{n-1} \sum_{1 \leq k < \ell \leq n} \left( \text{Ric}_{k \overline{k}}^{(1)} + \text{Ric}_{\ell \overline{\ell}}^{(1)} \right) +  \frac{1}{n-1} \sum_{1 \leq k < \ell \leq n} \left( \text{Ric}_{k \overline{k}}^{(2)} + \text{Ric}_{\ell \overline{\ell}}^{(2)} \right) \\
&=& \frac{1}{n-1} \sum_{1 \leq k < \ell \leq n} \left( \text{Ric}_{k \overline{k}}^{(1)} + \text{Ric}_{\ell \overline{\ell}}^{(1)} + \text{Ric}_{k \overline{k}}^{(2)} + \text{Ric}_{\ell \overline{\ell}}^{(2)} \right) \\
& \geq & \frac{2}{n-1} \sum_{1 \leq k < \ell \leq n} \left( \text{R}_{k \overline{\ell} \ell \overline{k}} + \text{R}_{\ell \overline{k} k \overline{\ell}} \right).
\end{eqnarray*}
\end{proof}

We consider the altered variant of the QOBC. 

\subsection*{Definition 3.7}
Let $(X,\text{g})$ be a Hermitian manifold. The \textit{altered quadratic orthogonal bisectional curvature} is defined \begin{eqnarray}
\widetilde{\text{QOBC}}_{\text{g}} :  \mathcal{F}_X \times \mathbf{R}^n_+ \to \mathbf{R},  \qquad \text{QOBC}_{\text{g}}(\xi) \ = \ \frac{1}{| \xi |_{\text{g}}^2} \sum_{\alpha, \gamma=1}^n \text{R}_{\alpha \bar{\gamma} \gamma \bar{\alpha}}(\xi_{\alpha} - \xi_{\gamma})^2.
\end{eqnarray} We have the following analog of \nameref{36} for the altered quadratic orthogonal bisectional curvature. 

\subsection*{Proposition 3.8}\label{37}
Let $(X,\text{g})$ be a Hermitian manifold with $\widetilde{\text{QOBC}}_{\text{g}} \geq 0$.  Then for any unitary pair of $(1,0)$-tangent vectors $\xi, \nu$,  we have \begin{eqnarray*}
\text{Ric}_{\text{g}}^{(3)}(\xi, \bar{\xi}) + \text{Ric}_{\text{g}}^{(3)}(\nu, \bar{\nu}) + \text{Ric}_{\text{g}}^{(4)}(\xi, \bar{\xi}) + \text{Ric}_{\text{g}}^{(4)}(\nu, \bar{\nu}) & \geq & 2 \left(\text{R}(\xi, \bar{\nu}, \nu, \bar{\xi}) + \text{R}(\nu, \bar{\xi}, \xi, \bar{\nu}) \right) .
\end{eqnarray*} In particular, in any local unitary frame,  the altered scalar curvature satisfies \begin{eqnarray*}
\widetilde{\text{Scal}}_{\text{g}} & \geq & \frac{1}{n-1} \sum_{1 \leq k < \ell \leq n} \left( \text{R}_{k \bar{k} \ell \bar{\ell}} + \text{R}_{\ell \bar{\ell} k \bar{k}} \right).
\end{eqnarray*} \begin{proof}
Suppose $\widetilde{\text{QOBC}}_{\text{g}} \geq 0$. Then in each unitary frame,  $\sum_{i,k =1}^n \text{R}_{i \overline{k} k \overline{i}}(\xi_i - \xi_j)^2 \ \geq \ 0$, for all $\xi = (\xi_1, ..., \xi_n) \in \mathbf{R}^n$.  Let $\xi_k =0$, $\xi_{\ell} =2$, and $\xi_j=1$ for $k \neq \ell$, $j \neq k$. Then \begin{eqnarray*}
4\left( \text{R}_{k \overline{\ell} \ell \overline{k}} + \text{R}_{\ell \overline{k} k \overline{\ell}} \right) + \sum_{j \neq k, j \neq \ell} \left( \text{R}_{k \overline{j} j \overline{k}} + \text{R}_{j \overline{k} k \overline{j}} + \text{R}_{\ell \overline{j} j \overline{\ell}} + \text{R}_{j \overline{\ell} \ell \overline{j}} \right) & \geq & 0.
\end{eqnarray*}
Let $f_k = \frac{1}{\sqrt{2}}(e_k - e_{\ell})$, $f_{\ell} = \frac{1}{\sqrt{2}}(e_k + e_{\ell})$, and $f_j = e_j$. Then \begin{eqnarray*}
&& \text{R}(e_k - e_{\ell}, \overline{e_k + e_{\ell}}, e_k + e_{\ell}, \overline{e_k - e_{\ell}}) + \text{R}(e_k + e_{\ell}, \overline{e_k - e_{\ell}}, e_k - e_{\ell}, \overline{e_k + e_{\ell}}) \\
&& + \frac{1}{2} \sum_{j \neq k, j \neq \ell} \left( \text{R}(e_k - e_{\ell}, \overline{e_j}, e_j, \overline{e_k - e_{\ell}}) + \text{R}(e_j, \overline{e_k - e_{\ell}}, e_k - e_{\ell}, \overline{e_j}) \right) \\
&& + \frac{1}{2} \sum_{j \neq k, j \neq \ell} \left( \text{R}(e_k + e_{\ell}, \overline{e_j}, e_j, \overline{e_k + e_{\ell}}) + \text{R}(e_j, \overline{e_k + e_{\ell}}, e_k + e_{\ell}, \overline{e_j}) \right) \\
&=& \text{R}_{k \overline{k} k \overline{k}} + \text{R}_{\ell \overline{\ell} \ell \overline{\ell}} - \text{R}_{k \overline{k} \ell \overline{\ell}} - \text{R}_{k \overline{\ell} k \overline{\ell}} + \text{R}_{k \overline{\ell} \ell \overline{k}} - \text{R}_{\ell \overline{k} \ell \overline{k}} + \text{R}_{\ell \overline{k} k \overline{\ell}} - \text{R}_{\ell \overline{\ell} k \overline{k}} \\
&& - \text{R}_{k \overline{k} k \overline{\ell}} + \text{R}_{k \overline{k} \ell \overline{k}} + \text{R}_{k \overline{\ell} k \overline{k}} - \text{R}_{\ell \overline{k} k \overline{k}} - \text{R}_{k \overline{\ell} \ell \overline{\ell}} + \text{R}_{\ell \overline{k} \ell \overline{\ell}} + \text{R}_{\ell \overline{\ell} k \overline{\ell}} - \text{R}_{\ell \overline{\ell} \ell \overline{k}} \\
&& - \text{R}_{k \overline{k} k \overline{\ell}} + \text{R}_{k \overline{k} \ell \overline{k}} - \text{R}_{k \overline{\ell} k \overline{k}} + \text{R}_{\ell \overline{k} k \overline{k}} + \text{R}_{k \overline{\ell} \ell \overline{\ell}} - \text{R}_{\ell \overline{k} \ell \overline{\ell}} - \text{R}_{\ell \overline{\ell} k \overline{\ell}} + \text{R}_{\ell \overline{\ell} \ell \overline{k}} \\
&& + \text{R}_{k \overline{k} k \overline{k}} + \text{R}_{\ell \overline{\ell} \ell \overline{\ell}} - \text{R}_{k \overline{k} \ell \overline{\ell}} - \text{R}_{k \overline{\ell} k \overline{\ell}} + \text{R}_{k \overline{\ell} \ell \overline{k}} - \text{R}_{\ell \overline{k} \ell \overline{k}} + \text{R}_{\ell \overline{k} k \overline{\ell}} - \text{R}_{\ell \overline{\ell} k \overline{k}} \\
&& + \text{R}_{k \overline{k} k \overline{\ell}} - \text{R}_{k \overline{k} \ell \overline{k}} - \text{R}_{k \overline{\ell} k \overline{k}} + \text{R}_{\ell \overline{k} k \overline{k}} + \text{R}_{k \overline{\ell} \ell \overline{\ell}} - \text{R}_{\ell \overline{k} \ell \overline{\ell}} - \text{R}_{\ell \overline{\ell} k \overline{\ell}} + \text{R}_{\ell \overline{\ell} \ell \overline{k}} \\
&& + \frac{1}{2} \sum_{j \neq k, j \neq \ell} \left(  \text{R}_{k \overline{j} j \overline{k}} + \text{R}_{\ell \overline{j} j \overline{\ell}} - \text{R}_{k \overline{j} j \overline{\ell}} - \text{R}_{\ell \overline{j} j \overline{k}} + \text{R}_{j \overline{k} k \overline{j}} + \text{R}_{j \overline{\ell} \ell \overline{j}} - \text{R}_{j \overline{k} \ell \overline{j}} - \text{R}_{j \overline{\ell} k \overline{j}}    \right) \\
&& + \frac{1}{2} \sum_{j \neq k, j \neq \ell} \left(  \text{R}_{k \overline{j} j \overline{k}} + \text{R}_{\ell \overline{j} j \overline{\ell}} + \text{R}_{k \overline{j} j \overline{\ell}} + \text{R}_{\ell \overline{j} j \overline{k}} + \text{R}_{j \overline{k} k\overline{j}} + \text{R}_{j \overline{\ell} \ell \overline{j}} + \text{R}_{j \overline{k} \ell \overline{j}} + \text{R}_{j \overline{\ell} k \overline{j}}     \right) \\
&=& 2 \left(  \text{R}_{k \overline{k} k \overline{k}} + \text{R}_{\ell \overline{\ell} \ell \overline{\ell}} - \text{R}_{k \overline{k} \ell \overline{\ell}} - \text{R}_{k \overline{\ell} k \overline{\ell}} + \text{R}_{k \overline{\ell} \ell \overline{k}} - \text{R}_{\ell \overline{k} \ell \overline{k}} + \text{R}_{\ell \overline{k} k \overline{\ell}} - \text{R}_{\ell \overline{\ell} k \overline{k}}   \right)  \\
&& + \sum_{j \neq k, j \neq \ell} \left( \text{R}_{k \overline{j} j \overline{k}} + \text{R}_{\ell \overline{j} j \overline{\ell}} + \text{R}_{j \overline{k} k \overline{j}} + \text{R}_{j \overline{\ell} \ell \overline{j}}  \right)  \ \geq \ 0.
\end{eqnarray*} Similarly, setting $f_k = \frac{1}{\sqrt{2}}(e_k - \sqrt{-1} e_{\ell})$, $f_{\ell} = \frac{1}{\sqrt{2}}(e_k + \sqrt{-1} e_{\ell})$, $f_j = e_j$, we have \begin{eqnarray*}
&& \text{R}(e_k - \sqrt{-1} e_{\ell}, \overline{e_k + \sqrt{-1} e_{\ell}}, e_k + \sqrt{-1} e_{\ell}, \overline{e_k - \sqrt{-1} e_{\ell}}) \\
&& + \text{R}(e_k + \sqrt{-1} e_{\ell}, \overline{e_k - \sqrt{-1} e_{\ell}}, e_k - \sqrt{-1} e_{\ell}, \overline{e_k + \sqrt{-1} e_{\ell}}) \\
&& + \frac{1}{2} \sum_{j \neq k, j \neq \ell} \left( \text{R}(e_k - \sqrt{-1} e_{\ell}, \overline{e_j}, e_j, \overline{e_k - \sqrt{-1} e_{\ell}}) + \text{R}(e_j, \overline{e_k - \sqrt{-1} e_{\ell}}, e_k - \sqrt{-1} e_{\ell}, \overline{e_j}) \right) \\
&& + \frac{1}{2} \sum_{j \neq k, j \neq \ell} \left( \text{R}(e_k + \sqrt{-1} e_{\ell}, \overline{e_j}, e_j, \overline{e_k + \sqrt{-1} e_{\ell}}) + \text{R}(e_j, \overline{e_k + \sqrt{-1} e_{\ell}}, e_k + \sqrt{-1} e_{\ell}, \overline{e_j}) \right) \\
&=& 2 \left( \text{R}_{k \overline{k} k \overline{k}} + \text{R}_{\ell \overline{\ell} \ell \overline{\ell}} - \text{R}_{k \overline{k} \ell \overline{\ell}} + \text{R}_{k \overline{\ell} k \overline{\ell}} + \text{R}_{k \overline{\ell} \ell \overline{k}} + \text{R}_{\ell \overline{k} \ell \overline{k}} + \text{R}_{\ell \overline{k} k \overline{\ell}} - \text{R}_{\ell \overline{\ell} k \overline{k}} \right) \\
&& + \sum_{j \neq k, j \neq \ell} \left( \text{R}_{k \overline{j} j \overline{k}} + \text{R}_{\ell \overline{j} j \overline{\ell}} + \text{R}_{j \overline{k} k \overline{j}} + \text{R}_{j \overline{\ell} \ell \overline{j}} \right) \ \geq \ 0.
\end{eqnarray*}
Hence, we see that \begin{eqnarray*}
&& 2 \left( \text{R}_{k \overline{k} k \overline{k}} + \text{R}_{\ell \overline{\ell} \ell \overline{\ell}} - \text{R}_{k \overline{k} \ell \overline{\ell}} - \text{R}_{\ell \overline{\ell} k \overline{k}} + \text{R}_{k \overline{\ell} \ell \overline{k}} + \text{R}_{\ell \overline{k} k \overline{\ell}} \right) \\
&& \hspace*{3cm} + \sum_{j \neq k, j \neq \ell} \left( \text{R}_{k \overline{j} j \overline{k}} + \text{R}_{\ell \overline{j} j \overline{\ell}} + \text{R}_{j \overline{k} k \overline{j}} + \text{R}_{j \overline{\ell} \ell \overline{j}} \right) \ \geq \ 0.
\end{eqnarray*}
Since \begin{eqnarray*}
&& \text{Ric}_{k \overline{k}}^{(3)} + \text{Ric}_{\ell \overline{\ell}}^{(3)} + \text{Ric}_{k \overline{k}}^{(4)} + \text{Ric}_{\ell \overline{\ell}}^{(4)} \\
&=& \text{R}_{k \overline{k} k \overline{k}} + \text{R}_{\ell \overline{\ell} \ell \overline{\ell}} + \text{R}_{\ell \overline{k} k \overline{\ell}} + \text{R}_{k \overline{\ell} \ell \overline{k}} + \sum_{j \neq k, j \neq \ell} (\text{R}_{j \overline{k} k \overline{j}} + \text{R}_{j \overline{\ell} \ell \overline{j}}) \\
&& + \text{R}_{k \overline{k} k \overline{k}} + \text{R}_{\ell \overline{\ell} \ell \overline{\ell}} + \text{R}_{k \overline{\ell} \ell \overline{k}} + \text{R}_{\ell \overline{k} k \overline{\ell}} + \sum_{j \neq k, j \neq \ell} (\text{R}_{k \overline{j} j \overline{k}} + \text{R}_{\ell \overline{j} j \overline{\ell}}) \\
&=& 2 \left( \text{R}_{k \overline{k} k \overline{k}} + \text{R}_{\ell \overline{\ell} \ell \overline{\ell}} + \text{R}_{\ell \overline{k} k \overline{\ell}} + \text{R}_{k \overline{\ell} \ell \overline{k}} \right) + \sum_{j \neq k, j \neq \ell} (\text{R}_{j \overline{k} k \overline{j}} + \text{R}_{j \overline{\ell} \ell \overline{j}} + \text{R}_{k \overline{j} j \overline{k}} + \text{R}_{\ell \overline{j} j \overline{\ell}} ),
\end{eqnarray*} it follows that \begin{eqnarray*}
\text{Ric}_{k \overline{k}}^{(3)} + \text{Ric}_{\ell \overline{\ell}}^{(3)} + \text{Ric}_{k \overline{k}}^{(4)} + \text{Ric}_{\ell \overline{\ell}}^{(4)} & \geq & 2 \left( \text{R}_{k \overline{k} \ell \overline{\ell}} + \text{R}_{\ell \overline{\ell} k \overline{k}} \right).
\end{eqnarray*} For the statement concerning the altered scalar curvature, simply observe that \begin{eqnarray*}
2 \widetilde{\text{Scal}}_{\text{g}} &=& \frac{1}{n-1} \sum_{k < \ell} (\text{Ric}_{k \overline{k}}^{(3)} + \text{Ric}_{\ell \overline{\ell}}^{(3)}) + \frac{1}{n-1} \sum_{k< \ell} (\text{Ric}_{k \overline{k}}^{(4)} + \text{Ric}_{\ell \overline{\ell}}^{(4)})\\
&=& \frac{1}{n-1}  \sum_{1 \leq k < \ell \leq n} \left( \text{Ric}_{k \overline{k}}^{(3)} + \text{Ric}_{\ell \overline{\ell}}^{(3)} + \text{Ric}_{k \overline{k}}^{(4)} + \text{Ric}_{\ell \overline{\ell}}^{(4)} \right) \\
& \geq & \frac{2}{n-1} \sum_{1 \leq k < \ell \leq n} (\text{R}_{k \overline{k} \ell \overline{\ell}} + \text{R}_{\ell \overline{\ell} k \overline{k}}).
\end{eqnarray*}
\end{proof}

\noindent We have the following immediate corollary. 

\subsection*{Corollary 3.9}\label{38}
Let $(X,\text{g})$ be a Hermitian manifold.  If $\text{QOBC}_{\text{g}} \geq 0$ and $\widetilde{\text{QOBC}}_{\text{g}} \geq 0$,  then $\text{Scal}_{\text{g}} \geq 0$ and $\widetilde{\text{Scal}}_{\text{g}} \geq 0$.  If,  in addition,  $\widetilde{\text{QOBC}}_{\text{g}} >0$,  then $\text{Scal}_{\text{g}} >0$ and the Kodaira dimension is negative.

\end{document}